\documentclass[12pt]{article}

\usepackage{eucal}
\usepackage{amsfonts,amssymb,amsmath}

\usepackage[latin1]{inputenc}
\usepackage{amsthm}
\usepackage{a4wide}
\usepackage{epsfig,graphics}
\usepackage{hyperref}
\usepackage{color}

\newcommand{\x}{{{\mathbf x}}}
\newcommand{\e}{{{\mathbf e}}}
\newcommand{\m}{{{\mathbf m}}}
\newcommand{\tm}{{ \widetilde{\mathbf m}}}
\newcommand{\tmu}{{ \widetilde{\mu}}}
\newcommand{\dd}{{{\mathbf d}}}
\newcommand{\y}{{{\mathbf y}}}
\newcommand{\rr}{{{\mathbf r}}}
\newcommand{\oo}{{{\mathbf o}}}

\newcommand{\XX}{{{\mathbf X}}}
\newcommand{\z}{{{\mathbf z}}}
\newcommand{\vb}{{{\mathbf v}}}
\newcommand{\aA}{{{A}}}
\newcommand{\aB}{{{B}}}
\newcommand{\aC}{{{C}}}
\newcommand{\aD}{{{D}}}

\newcommand{\Z}{{{\mathbb{Z}^{d}}}}

\newcommand{\R}{{{\mathbb{R}}}}
\newcommand{\ed}{\mathrm{d}}
\newcommand{\kK}{\mathrm{K}}

\newcommand{\cN}{{{\mathcal{N}}}}
\newcommand{\cF}{{{\mathcal{F}}}}

\newcommand{\E}{\text{\sf E}}
\renewcommand{\Pr}{\text{\sf P}}
\newcommand{\ovl}[1]{\overline{#1}}
\newcommand{\reff}[1]{(\ref{#1})}

\newcommand{\one}{\mathbf{1}}

\newtheorem{theorem}{Theorem}[section]
\newtheorem{lemma}[theorem]{Lemma}

\newtheorem{iten}[theorem]{}

\begin{document}

\title{Harness Processes and Non-Homogeneous Crystals}

\author{Pablo A. Ferrari, Beat M. Niederhauser, Eugene A. Pechersky}

\date{}
\maketitle

\begin{abstract}
  We consider the Harmonic crystal, a measure on $\mathbb{R}^{\mathbb{Z}^{d}}$
  with Hamiltonian $H(\x)=\sum_{i,j}J_{i,j}(\x(i)-\x(j))^{2}+
  h\sum_{i}(\x(i)-\dd(i))^{2}$, where $\x,\,\dd$ are configurations,
  $\x(i),\dd(i)\in\mathbb{R}$, $i,j\in{\mathbb{Z}^{d}}$. The configuration $\dd$
  is given and considered as observations. The `couplings' $J_{i,j}$ are finite
  range. We use a version of the harness process to explicitly construct the
  unique infinite volume measure at finite temperature and to find the unique
  ground state configuration $\m$ corresponding to the Hamiltonian.
\end{abstract}


\section{Introduction}
The harnesses were introduced by Hammersley \cite{H} to model the
behavior of a crystal and to introduce a multi-dimension version of a
martingale.  Let $P=(p(i,j),\,i,j\in\Z)$ be a homogeneous symmetric
Markov transition matrix with
$p(i,i)=0$. A \emph{harness} is a measure on
$\mathbb{R}^{\mathbb{Z}^{d}}$ with the property that the mean height
at $i$ given the heights at all sites different of $i$ is a
$P$-weighted mean of the heights of the other sites. The \emph{serial
  harness} is a Markov process on $\mathbb{R}^{\mathbb{Z}^{d}}$
updated at all discrete times at all sites by the rule: substitute the
height at site $i$ by a $P$-weighted mean of the neighbors plus a
centered independent random variable (the noise). Hsiao \cite{Hs1}
proposed a continuous-time version then called \emph{harness process}
in \cite{FN}. The heights are updated at Poisson epochs using the same
rule as in the serial harness. If the noise is a centered Gaussian
random variable, the reversible measure of the process is the harmonic
crystal, that is, the Gibbs measure with Hamiltonian
\begin{equation}
  \label{p6}
H(\x):=\sum_{i,j}J_{i,j}(\x(i)-\x(j))^{2}
\end{equation}
where $\x\in\R^{\Z}$, $\x(i)$ for $i\in{\Z}$ represents the height at site $i$
and $J_{i,j} = p(i,j)$; the temperature $1/\beta$ is given by the variance of
the noise.

We study a version of the harness process with a external local data term.  We
can think that each site $i\in\mathbb{Z}^{d}$ has an additional ``neighbor''
with a fixed height $\dd(i)$, the data. The updating of the height at $i$
involves the data $\dd(i)$ in the averaging. This is a ``heat bath'' dynamics
associated to the quadratic (Gaussian) Hamiltonian \eqref{16}: at rate 1, the
height at site $i$ is substituted by a random height distributed with the
conditional distribution associated to the Hamiltonian, given the heights at the
other sites. We show ergodicity of the harness process with data $\dd\in\XX$, a
set of configurations with a mild restriction on the growth, defined in
\eqref{e13}; this extends the work of Hsiao \cite{Hs1, Hs2} who considered the
case $\dd(i)\equiv 0$. Ergodicity means that there exists a unique invariant
measure for the process and that the process starting from measures
concentrating mass on $\XX$ converges to the invariant measure.  The unique
invariant measure is also reversible. We show that any infinite-volume Gibbs
measure on $\XX$ associated to the quadratic Hamiltonian \eqref{16} is invariant
for the dynamics. This fact and ergodicity imply that there is only one
infinite-volume Gibbs measure associated to the Hamiltonian. When the process is
constructed in a finite subset $\Lambda$ of $\Z$, the invariant measure is a
harmonic crystal with external site-by-site field $\dd$. As $\Lambda$ grows to
$\Z$ the harmonic crystals in $\Lambda$ converge to the unique invariant measure
for the harness process in $\Z$. If the data is flat (i.e.  constant) then we
are in the case of massive lattice models in quantum field theory which are well
known (see \cite{GJ, GJS, D}).

Shortly our method is as following. We slice the space-time configuration space
$\Xi=\R^{\Z}\times\R_+$ in pieces determined by the realizations of the Poisson
processes governing the updating times. We show the convergence of the process
almost surely for (almost) every slice. To this end we use Harris graphical
construction of the process as a function of a space-time marked Poisson process
of rate 1 on $\Z\times\R$. It is convenient to construct the process in an
arbitrary interval $[s,t]$ to be able to take the limit as $s\to-\infty$ which
is equivalent to take the limit as $t\to\infty$ in distribution. At each Poisson
epoch, the value of the process at corresponding site is substituted by an
average of the values of the process at the neighboring sites plus an
independent noise. When the construction is explored backwards in time the value
of the process at site $i$ at time $t$ can be expressed in function of the
probabilities of a random walk running backwards in time \emph{conditioned} on
the space-time epochs of the Poisson process. The walk is killed at a rate
related to the weight of the external data and when it hits the boundary. The
value of the process at site $i$ at time $t$ starting at time $s$ is expressed
as a sum of four terms: the contributions given by (1) the noise, (2) the
external data, (3) the boundary condition (in case the process is studied in a
finite region) and (4) the initial condition.  The noise contribution is a
martingale with uniformly bounded second moments so it converges as $s\to
-\infty$. The data contribution converges to a deterministic function $\m$, a
harmonic function for a kernel associated to $J$ and $h$.  The boundary and
initial contributions go to zero as the region grows to $\Z$ and the initial
time $s$ goes to $-\infty$, respectively.

In Section \ref{S2} we state the main three theorems.  The first one shows the
existence of the Harness process in infinite volume.  The second one says that
the Harness process with external data is ergodic, that its unique invariant
measure coincides with the infinite-volume Gibbs measure and that there is a
unique infinite volume Gibbs measure for this Hamiltonian. The third theorem
shows that the harmonic function $\m$ is the unique function minimizing the
Hamiltonian with the data term (ground state).  In Section \ref{S3} we construct
the process and state intermediate results for the finite and
infinite process.  In Sections \ref{S4} and \ref{S5} we show that the finite
process has as reversible measure the Gibbs measure in finite volume and that
the space and time limits coincide.  In Section \ref{S6} we show that the
harmonic function $\m$ is the unique minimizer.

\section{Main Results}
\label{S2}
Let $\Lambda\subset\Z$ be a finite set and define the Hamiltonian
\begin{equation}\label{16}
  H_\Lambda(\x)=\sum_{i,j: \{i,j\}\not\subset \Lambda^c}J_{i,j}(\x(i)-\x(j))^{2}+
  h\sum_{i\in\Lambda}(\x(i)-\dd(i))^{2}
\end{equation}
$\x,\dd\in\mathbb{R}^{\mathbb{Z}^{d}}$, $J_{i,j}$ is a finite range pair
potential, that is $J_{i,j}=0$ if $|i-j|\geq R>0$, $\dd$ is a fixed
configuration which can be taken as data and $h>0$ is a fixed parameter.  In
fact $H_\Lambda=H_\Lambda(\cdot,\dd,h)$ depends on $\dd$ and $h$ which are fixed
and dropped from the notation unless necessary. Let $\x_{\Lambda}$ be the
configuration $\x$ restricted to the set $\Lambda$ and the superposition
configuration $\x_\Lambda\y_{\Lambda^c}$ be given by $\x$ for sites in $\Lambda$
and by $\y$ otherwise. Let $|\Lambda|$ be the number of sites in $\Lambda$. The
family of measures $\{\mu_{\Lambda}(\cdot|\y_{\Lambda^{c}}),\,
\Lambda\subset\Z,\,|\Lambda|<\infty,\, \y\in\R^\Z\}$, defined by
\begin{equation}
\label{10}
\mu_{\Lambda}(\ed\x_{\Lambda}|\y_{\Lambda^{c}})= \frac{
  e^{-H_\Lambda(\x_\Lambda\y_{\Lambda^c})}}{Z_{\Lambda}}\mu_{\Lambda}^{0}(\ed\x_{\Lambda}),
\end{equation}
where $\mu_{\Lambda}^{0}$ is the Lebesgue measure in $\mathbb{R}^{\Lambda}$ and
$Z_{\Lambda}$ is a normalizing constant, is called the specification
associated to the Hamiltonian $H$. A measure $\mu$ defined on $\R^\Z$ is said to
satisfy the Dobrushin-Lanford-Ruelle (DLR) equations if its conditioned
distributions coincide with the specifications:
\begin{equation}
  \label{dlr1}
  \mu(\ed\x_{\Lambda}|\x_{\Lambda^{c}}=\y_{\Lambda^{c}})=\mu_{\Lambda}(\ed\x_{\Lambda}|\y_{\Lambda^{c}})
\end{equation}
for $\mu$-almost all $\y$.

We construct the dynamics in a set of configurations with limited growth. Let
\begin{equation}
  \label{e13}
  \XX = \left\{\x \in \R^{\Z}\,:\, \sum_j|\x(j)|\alpha^{|j|/R} <\infty
    \hbox{ for all }i\in\Z\right\}
\end{equation}
Recall that $R$ is the radius of the interactions.  We introduce the {\it
  harness process} $\{\eta_{t}(i),\;i\in \mathbb{Z}^{d}\}$ on $\XX$. Let
$\alpha\in[0,1]$, $P=(p(i,j),\,i,j\in\mathbb{Z}^{d})$ be a space homogeneous
finite-range symmetric stochastic matrix (that is, $p(i,j)\geq
0,\;\sum_{j}p(i,j)=1$ for all $i$, $p(i,i+j)=p(0,j),\; p(0,j)=0 \mbox{ if
}|j|\geq R>0$ and $p(i,j)=p(j,i)$). Symmetry is not necessary to define the
process but it is natural in this context: we relate later $p(i,j)$ with
$J_{i,j}$ which is symmetric.

The generator of the process acts on locally finite continuous functions $f$ by
\begin{equation}
\label{pb}
Lf(\x)=\sum\limits_{k\in\mathbb{Z}^{d}}L_{k}f(\x)=
\sum\limits_{k\in\mathbb{Z}^{d}} \int G(\ed
x)\left[f\left(\alpha P_{k}\x+(1-\alpha)\e_{k}\dd(k)
  +\e_{k}x\right)-f(\x)\right],
\end{equation}
where $G(\ed x)= \frac{1}{\sqrt{2\pi}}e^{-x^{2}/2}\ed x$ is the standard
Gaussian distribution (with zero mean and variance 1), called \emph{noise} and
$P_{k}\x$ is the configuration defined by
\begin{equation}
\label{pb1}
(P_{k}\x)(j)=\begin{cases}
  \sum_{i\in\mathbb{Z}^{d}}p(k,i)\x(i)&\mbox{ if }j=k,\\
  \x(j)&\mbox{ otherwise }
\end{cases}
\end{equation}
and $\e_{k}(j)=\one\{k=j\}$. In this process at rate 1 the height at site $i$ is
updated with a convex combination of the heights at the neighbors of $i$ and the
height $\dd(i)$, plus an independent standard Gaussian variable. That is, if
site $i$ is updated at time $t$, the height at $i$ is substituted by a random
variable with the same law as
\begin{equation}
  \label{rf8}
  \alpha\, (P_i\eta_{t-})(i)+(1-\alpha)\dd(i) +Z
\end{equation}
where $Z$ is a standard Gaussian independent of $\eta_{t-}$. The real variable
\eqref{rf8} has law $\mu_{\{i\}}(\cdot\,|(\eta_{t-})_{\{i\}^c})$, the
distribution of the $\{i\}$ coordinate given the values of $\eta_{t-}$ in
$\{i\}^c$ as defined in \reff{10}.  Our dynamics coincides with the so called
``heat bath''.

The following three theorems are the aim of this work. The first one shows the
existence of the process.
\begin{theorem}
  \label{p3} Assume $\dd\in\XX$. There exists a Markov process $(\eta_t)$ on
  $\XX$ with generator~$L$:
\begin{equation}
\label{p4}
\lim_{u\rightarrow 0}\frac{1}{u}
\E[f(\eta_{t+u})-f(\eta_t)|\EuScript{F}_{t}]=Lf(\eta_t),
\end{equation}
where $\EuScript{F}_{t}$ is the sigma algebra generated by $(\eta_s,\, s\le t)$,
the past of $\eta_s$ up to time $t$.
\end{theorem}
Theorem \ref{p3} is a consequence of a general existence result of
Basis \cite{Ba1}; we provide here an alternative construction.

Let $S(t)$ be the semigroup associated to $L$ defined by $(S(t)f)(\eta) = \E
(f(\eta_t)|\eta_0=\eta)$. It acts on measures by $\int f\,\ed (\nu S(t)) =\int
S(t)f\, \ed \nu $.  Our second result says that the harmonic crystal is the
invariant measure for the harness process when $p(0,0)=0$ and establishes time
and space limits.

\begin{theorem}\label{tp6}
  Assume $\dd\in\XX$ and $p(0,0)=0$. (i) The following time and space limits
  exist and are identical.  For any initial measure $\nu$ concentrating on
  $\XX$, any boundary conditions $\y\in\XX$ and any increasing sequence
  $\Lambda\nearrow\mathbb{Z}^{d}$,
\begin{equation}
  \label{12} \lim_{t\to\infty}\nu S(t)
  =\lim_{\Lambda\nearrow\mathbb{Z}^{d}}\mu_{\Lambda}:=\mu.
\end{equation}
(ii) $\mu$ is reversible for $L_{k}$ for all $k\in\Z$, in particular it is
reversible for $L$. (iii) $\mu$ is the unique measure in $\XX$ satisfying the
DLR equations for the specifications (\ref{10}) with $J_{i,j}=\alpha p(i,j)$ and
$h=1-\alpha$.
\end{theorem}

The theorem implies that the harness process is ergodic in $\XX$: there exists a
unique invariant measure $\mu$ and the process starting in $\XX$ converges to
$\mu$. It also implies that for any boundary conditions in $\XX$ the
thermodynamic limit is unique (absence of phase transition).

The time convergence was proven by Hsiao \cite{Hs1} and \cite{Hs2} in the case
$\dd(k)\equiv 0$. The space convergence is contained in Spitzer \cite{Sp} and
Dobrushin \cite{D2}, see also Caputo \cite{C}.  Our approach permits to
construct simultaneously (coupling) realizations of the measures in all finite
boxes and the infinite volume measure in such a way that the convergence is
almost sure. Notice however that the uniqueness result is a consequence of the
fact that there is only one infinite-volume invariant measure for the process
and that any measure compatible with the DLR conditions is invariant. The space
thermodynamical limit is not used in the proof of
uniqueness. 

The condition $p(0,0)=0$ is necessary to guarantee that the measure
$\mu_\Lambda$ is invariant for the generator $L_k$. 

In Statistical Mechanics it is natural to extend the specifications \eqref{10}
to a family of measures $\mu^\beta_{\Lambda}$ defined by
\begin{equation}
  \label{m1}
  \mu^\beta_{\Lambda}(\ed\x_{\Lambda}|\y_{\Lambda^{c}})=
  \frac{1}{Z^\beta_{\Lambda}}\,
  e^{-\beta H_\Lambda(\x_\Lambda\y_{\Lambda^c})}\,\mu_{\Lambda}^{0}(\ed\x_{\Lambda}),
\end{equation}
for $\beta> 0$; $\beta$ is called the inverse temperature. We consider in detail
the case $\beta=1$; the other cases reduce to this one using $\beta
H_\Lambda(\x_\Lambda\y_{\Lambda^c},\dd,h)=
H_\Lambda(\x^\beta_\Lambda\y^\beta_{\Lambda^c},\dd^\beta,h)$, where $\x^\beta(i)
= \sqrt\beta\x(i)$, etc.

When $\beta=\infty$ the randomness vanishes and $\mu^\infty_\Lambda$ is
interpreted as the measure concentrating mass on configurations $\x_\Lambda$
minimizing $H_{\Lambda}(\x_\Lambda\y_{\Lambda^c})$ for finite $\Lambda$. When
$\Lambda=\Z$, we denote $H(\x)=H_{\Z}(\x)$, an infinite sum only formally
defined. In this case we need to give a sense to the word ``minimizing''. If
$\tilde\x$ differs from $\x$ on a finite set of sites $\Lambda$, then the
infinite sums defining $H(\x)$ and $H(\tilde\x)$ differ only on a finite number
of summands. We define $H(\tilde\x)-H(\x)$ as the difference of the
corresponding different summands, that is,
$H_\Lambda(\tilde\x_\Lambda\x_{\Lambda^c})-H_\Lambda(\x_\Lambda\x_{\Lambda^c})$.
We say that $\x$ \emph{minimizes} $H(\cdot)$ if $H(\tilde\x)-H(\x)>0$ for all
$\tilde\x$ local modification of $\x$. Measures defined on $\XX$ compatible with
the specifications \reff{10} with $\beta=\infty$ are called \emph{ground
  states}; see appendix B of van Enter, Fern\'andez and Sokal \cite{EFS} for
details. The ground states are concentrated on minimizing configurations.

Let $\kK(i,j)$ be the probability that the walk with rates $\alpha P$ (a walk
killed at rate $1-\alpha$) is killed at site $j$ when starting at site~$i$.  In
the next theorem we show that a delta measure concentrating mass in the
configuration given by the $\kK$-average of $\dd$ is the unique ground state.

\begin{theorem}\label{14}
  Assume $p(0,0)=0$ and that $\dd\in\XX$ and let $\m$ be the configuration given
  by
\begin{equation}
\label{15} \m(i):=\sum_j \kK(i,j)\dd(j) \qquad \hbox{ for all
} i\in \Z
\end{equation}
then $\m$ minimizes the Hamiltonian $H$ and the delta measure concentrating mass
on~$\m$ is the unique ground state for specifications \reff{10} in $\XX$.
\end{theorem}

Using the Kolmogorov Backwards equation for the walk with transitions $\alpha P$
killed at rate $1-\alpha$ and boundary conditions $\dd$, we see that $\m$
satisfies the equation
\begin{equation}
  \label{hh1}
  \m(i)\,=\, \sum_{j\in\Z}  \alpha p(i,j)\m(j) + (1-\alpha) \dd(i) .
\end{equation}
That is, $\m$ is a harmonic function for a transition matrix associated to $P$
and $\alpha$ in an extended graph with ``boundary conditions'' $\dd$. The
extended graph has vertices in $\Z\cup(\Z)^*$ where $(\Z)^*$ is a copy of $\Z$,
and edges $\big\{(i,j)\in\Z\times\Z\,:\,|i-j|=1\}\cup\{(i,i^*)\,:\, i\in \Z$ and
$i^*$ is the copy of $i$ in $(\Z)^*\big\}$. The boundary conditions are fixed in
$(\Z)^*$ equal to $\dd$; $\dd(i)$ is the value of the boundary condition at
$i^*$ and $(1-\alpha)$ is the weight of the edge $(i,i^*)$. $\alpha p(i,j)$ is
the weight of the edge $(i,j)$.

\section{Harris graphical construction}\label{S3}
The proof of the above theorems are based on an adaptation of the Harris
graphical construction for the harness process proposed by the first two authors
in \cite{FN}.  The core of the construction is a rate-one space-time Poisson
process $\cN$ on $\Z\times\R$.  This can be thought of as a product of
homogeneous one-dimensional Poisson process in $\R$, one for each $i\in\Z$.
Space-time points in $\cN$ are denoted $(i,\tau)$ and called \emph{epochs}. To
each $(i,\tau)\in\cN$ we attach two independent \emph{marks}: $\xi(i,\tau)$ and
$\varphi(i,\tau)$, where $\xi(i,\tau)$ is a Gaussian random variable with zero
mean and variance $1$ called \emph{noise} and $\varphi(i,\tau)$ are variables
whose distribution is described later; these two families are iid and mutually
independent and independent of $\cN$. We denote $\Pr$ and $\E$ the probability
and expectation induced by these Poisson processes with marks.

The Harness process is realized as a function of the Poisson epochs and the
marks $\xi$ (for the moment we do not use the marks $\varphi$) as follows.  The
height at each site $i$ only changes at times $(i,s)\in \cN$.  Assuming that the
configuration at time $s-$ is $\eta_{s-}$ and $(i,s)\in\cN$ is a Poisson epoch,
then the configuration at time $s$ at sites $k\ne i$ does not change
($\eta_s(k)=\eta_{s-}(k)$) and
\begin{equation}
  \label{f7}
  \eta_s(i) = \alpha (P_i \eta_{s-})(i)+(1-\alpha) \dd(i)+ \xi(i,s)
\end{equation}
where $P_i$ is defined in \eqref{pb1}.  In other words, at the Poisson epoch
$(i,s)$, the height at $i$ is substituted by a average of the other sites and
the external data at $i$ plus a Normal random variable independent of
``everything'' else. By a standard percolation argument this construction can be
performed for small time intervals, so that $\Z$ is partitioned in non
interacting pieces); we sketch it in the proof of \ref{i4a} in Section \ref{iv1}
later. We give another construction based on a ``dual'' representation using the
variables $\varphi$.

Let the variables $\varphi(i,\tau)$ be independent with law
\begin{align}\label{13.2}
\Pr(\varphi(i,\tau)=j)&=\alpha p(i,j),\\
\Pr(\varphi(i,\tau)=i^{*})&=1-\alpha,\nonumber
\end{align}
where $i^*$ ($\not\in \Z$) is the copy of $i$ in $(\Z)^*$. The parameter $\alpha
\in (0,1)$ is later chosen as in Theorem \ref{tp6}. Define a family of backward
random walks indexed by $(\gamma,t)$, the space-time starting point, as a
deterministic function of the Poisson epochs $\cN$ and the marks $\varphi$ (here
we do not use the marks $\xi$). Fix $t\in\R$ and $\gamma\in \Z\cup(\Z)^*$. For
each $s<t$ we define $\sigma^\gamma_{[s,t]}\in \Z\cup(\Z)^*$ as the position of
a random walk going backwards in time with initial position (at time $t$)
$\sigma^\gamma_{[t,t]}=\gamma$ and evolving with the following rules.  The walk
does not move between Poisson epochs and for $s<t$, if at time $s+$ the walk is
at site $i\in\Z$ and $(i,s)\in\cN$, then at time $s$ the walk jumps to the
position $\varphi(i,s)$: for $s<t$,
\begin{equation}
  \label{f1}
  \sigma^\gamma_{[s,t]}= \begin{cases}
    \varphi(i,s)&\mbox{if }\sigma^\gamma_{[s+,t]}=i\in\Z\hbox{ and }(i,s)\in \cN \\
    \sigma^\gamma_{[s+,t]}&\mbox{if }\sigma^\gamma_{[s+,t]}=i^*\in(\Z)^*
    \hbox{ or } \sigma^\gamma_{[s+,t]}=i\in\Z\hbox{ and } (i,s)\notin \cN
\end{cases}
\end{equation}
We say that the walk $\sigma^\gamma_{[s,t]}$ is \emph{absorbed} at
$j\in\mathbb{Z}^d$ at time $s$ if $\sigma^\gamma_{[u,t]}=j^*$ for $u\le s$ and
$\sigma^\gamma_{[u,t]}\in\Z$ for $u\in(s,t]$.  The family
$((\sigma^\gamma_{[u,t]},\, u\le t), \;t\in \R, \gamma\in \Z\cup(\Z)^*))$ is a
function of $\cN$ and $\varphi$, but we drop this dependence in the notation.

A key object in this analysis is the law of the walk $\sigma$ \emph{conditioned}
on a realization of the Poisson epochs $\cN$. For $i\in\Z$ let
\begin{equation}
  \label{f3}
 b_{[u,t]}(i,j)=\Pr\left(\sigma^i_{[u,t]}= j\big|\cN \right)
\end{equation}
be the probability the walk $\sigma$ starting at $i$ at time $t$ to be at $j$ at
time $u$ \emph{given} the Poisson epochs. The probabilities $b_{[u,t]}(i,j)$ are
function of the Poisson epochs $\cN$ and do \emph{not} depend on
$\varphi$. Analogously, for $(j,\tau)\in\cN$ define the probability of
absorption at time $\tau$ at site $j$ given the Poisson epochs by
\begin{equation}
  \label{f3a}
 a_{[\tau,t]}(i,j^*)=\Pr\left(\sigma^i_{[\tau,t]}= j^*,
 \sigma^i_{[u,t]}\in\Z,\hbox{ for }u\in(\tau,t]\big|\cN \right)
\end{equation}
Notice that
\begin{equation}
  \label{f4}
 a_{[\tau,t]}(i,j^*)= (1-\alpha)\,b_{[\tau,t]}(i,j)
\end{equation}
Define the $\XX$-valued process $\eta_{[s,t]}$ in the time interval $[s,t]$ with
$s\le t$ and initial condition $\z$ at time $s$ by
\begin{align}
  \nonumber
  \eta_{[s,s]}&\equiv\z,\\
  \label{3.1} \eta_{[s,t]}(i)&= \sum_{(j,\tau)\in \cN [s,t]}b_{[\tau,t]}(i,j)
  \xi(j,\tau)+a_{[\tau,t]}(i,j^*)\dd(j) + \sum_{j\in\Z}b_{[s,t]}(i,j) \z(j),
\end{align}
where $\cN [s,t]=\{(j,\tau)\in \cN:\, \tau\in[s,t]\}$ and recall $\xi(j,\tau)$
is the noise associated to the epoch $(j,\tau)$. Under this construction
$\eta_{[s,t]}(i)$ is a function of $\cN[s,t]$ and the corresponding noises
$\xi$; it is an average determined by $\cN$ of the noises $\xi$, the external
field $\dd$ and the initial condition $\z$ at time $s$.  Our goal is to prove
the following

\begin{iten}
  \label{i1} For each $s\in \R$ and $\z\in\XX$ the process $(\eta_{[s,t]},\,
  t\ge s)$ defined in (\ref{3.1}) is well-defined.  Namely, the sums in
  (\ref{3.1}) are finite with probability 1 and $\eta_{[s,t]}\in\XX$ for all
  $s<t$. Furthermore the process is Markovian with generator
  $L$ given in (\ref{pb}) and initial condition $\z$ at time $s$.
\end{iten}

\begin{iten}
  \label{i3} For any configuration $\z\in\XX $ and fixed $t\in \R$, the limit
  \begin{equation}
    \label{f6p}
    \lim_{s\to-\infty} \eta_{[s,t]}(i) :=\eta_t(i)
  \end{equation}
  exists with probability one and \emph{does not depend on $\z$}. The process
  $(\eta_t,\,t\in\R)$ is a stationary Markov process with generator $L$ given in
  (\ref{pb}).
\end{iten}

\begin{iten}
  \label{i4} Call $\mu$ the marginal law of $\eta_t$ (which does not depend on
  $t$), then $\mu$ satisfies the DLR equations \eqref{dlr1}.
\end{iten}

\begin{iten}
  \label{i4a} If $\tmu$ satisfies the DLR equations then $\tmu$ is
  reversible for the process $(\eta_{[s,t]},\,t\ge s)$.
\end{iten}

\begin{iten}
  \label{i4b} $\mu$, the marginal law of $\eta_t$, is the unique infinite volume
  Gibbs measure on $\XX$ for the specifications \eqref{10}.
\end{iten}

In other words, the strategy of our proof is to construct a stationary process
in infinite volume whose time marginal is the unique Gibbs measure for the
Hamiltonian \eqref{16}. Theorem \ref{p3} follows from \ref{i1}. 

\vspace{.5cm}

\paragraph{Finite volume}
We start considering the Hammersley processes in finite volume
$\Lambda\subset\Z$ with boundary conditions $\y\in\XX$.  The updates occurs at
space-time Poisson epochs $(i,s)\in\cN$ as follows:
\begin{equation}
  \label{87}
  \eta_s(i) = \alpha \Bigl(\sum_{j\in\Lambda}
  p(i,j)\eta_{s-}(j)+\sum_{j\in\Lambda^c}p(i,j)\y(j)\Bigr)+ (1-\alpha) \dd(i)+ \xi(i,s)
\end{equation}
That is, at the space-time Poisson epochs, the process substitutes the value at
$i$ by an average of the values at the other sites and the external data, the
values outside $\Lambda$ are kept fixed and given by the boundary configuration
$\y$.

The construction of the finite harness process in $\Lambda$ with boundary
configuration $\y$ goes along the same lines as in infinite volume, the
difference is that the probabilities $b$ and $a$ are computed for walks that are
also absorbed at $\Lambda^c$.

Call $\Lambda^*$ the copy of $\Lambda$ in $(\Z)^*$. We define a family of
backward random walks absorbed at $\Lambda^c\cup\Lambda^*$ indexed by
$(\gamma,[s,t],\Lambda)$, with space-time starting point $(\gamma,t)$ with
$\gamma\in\Lambda$. The rules now are
\begin{equation}
  \label{f1p}
  \sigma^\gamma_{[s,t],\Lambda}= \begin{cases}
    \varphi(i,s)&\mbox{if }\sigma^\gamma_{[s+,\,t],\Lambda}=i\in\Lambda
    \hbox{ and }(i,s)\in \cN \\
    \sigma^\gamma_{[s+,\,t],\Lambda}&\mbox{if }
    \sigma^\gamma_{[s+,\,t],\Lambda}\in(\Z)^* \cup \Lambda^c
    \hbox{ or } \sigma^\gamma_{[s+,\,t],\Lambda}=i\in\Z\hbox{ and } (i,s)\notin \cN
\end{cases}
\end{equation}
The only difference is that now the walk is absorbed at $\Lambda^*$ \emph{and} at
$\Lambda^c$. If the walk starts at $\Lambda$, it will be absorbed either at
$\Lambda^*$ or at $\Lambda^c$. Calling
\[
\cN([s,t],\Lambda) = \{(j,\tau)\in\cN:\, \tau\in[s,t],j\in\Lambda\},
\]
the family $((\sigma^\gamma_{[u,t],\Lambda},\, u\in[s, t]), \;t\in \R, \gamma\in
\Lambda$ is a function of $\cN([s,t],\Lambda)$ and the associated $\varphi$, but
we drop this dependence in the notation. For $i\in\Lambda$ and $u<t$ define
\begin{equation}
  \label{f3c}
 b_{[u,t],\Lambda}(i,j)=\Pr\left(\sigma^i_{[u,t],\Lambda}= j\big|\cN \right),
\end{equation}
the transition probabilities for the walk absorbed at $(\Z)^*$ and $\Lambda^c$
given the Poisson epochs. The probabilities of absorption at time $\tau$ at site
$\gamma\in\Lambda^*\cup\Lambda^c$ given the Poisson epochs are defined by
\begin{equation}
  \label{f3d}
 a_{[\tau,t],\Lambda}(i,\gamma)=\Pr\left(\sigma^i_{[\tau,t],\Lambda}= \gamma,
 \sigma^i_{[u,t],\Lambda}\in\Lambda,\hbox{ for }u\in(\tau,t]\big|\cN \right)
\end{equation}

We define the $\R^\Lambda$-valued process $\eta_{[s,t],\Lambda}$ in the time interval
$[s,t]$ with $s\le t$, initial condition $\z\in\XX$ and boundary
conditions $\y\in\XX$ at time $s$ by
\begin{align}
  \nonumber
  \eta_{[s,s],\Lambda}&\equiv\z_\Lambda \\
  \label{f31} \eta_{[s,t],\Lambda}(i) &= \sum_{(j,\tau)\in
    \cN([s,t],\Lambda)}\Bigl(b_{[\tau,t],\Lambda}(i,j)\xi(j,\tau) +
  a_{[\tau,t],\Lambda}(i,j^*)\dd(j)
  + \sum_{k\in \Lambda^c}a_{[\tau,t],\Lambda}(i,k)\,\y(k)\Bigr)\nonumber\\
  &\qquad\qquad+\; \sum_{j\in\Lambda}b_{[s,t],\Lambda}(i,j) \z_\Lambda(j),\quad \hbox{ for }
  i\in\Lambda
\end{align}
where $\cN ([s,t],\Lambda)=\{(j,\tau)\in \cN:\, j\in\Lambda,\, \tau\in[s,t]\}$,
the value $\xi(j,\tau)$ is the Gaussian random variable associated to the
Poisson epoch $(j,\tau)\in\cN$. 

With the above construction $\eta_{[s,t],\Lambda}(i)$ is an average determined
by $\cN$ (through the weights $a$ and $b$) of the noise $\xi$, the external
field $\dd$, the boundary configuration $\y$ and the initial condition $\z$ at
time $s$. We drop these dependences in the notation.

We prove the following facts about the finite process:

\begin{iten}
    \label{i5} For each $s\in\R$ the process $(\eta_{[s,t],\Lambda},\, t\ge
  s)$ is Markov with initial condition $\z_\Lambda$, boundary conditions
  $\y_{\Lambda^c}$ and generator
  \begin{equation}
    \label{p2f}
    \begin{aligned}
      L_{\Lambda}f(\x_\Lambda)&=\sum_{k\in\Lambda}L_{k}f(\x_\Lambda\y_{\Lambda^c})\\
      &=\;\sum_{k\in\Lambda} \int G(\ed x)\big[f\big([\alpha
      P_{k}(\x_\Lambda\y_{\Lambda^c})+(1-\alpha)\e_{k}\dd(k)
      +\e_{k}x]_\Lambda\big)-f(\x_\Lambda)\big],  
  \end{aligned}
\end{equation}
where $f$ is a bounded function depending only on coordinates in ${\Lambda}$
\end{iten}

\begin{iten}
   \label{i51}
   The measure $\mu_\Lambda(\cdot|\y_{\Lambda^c})$ given in \eqref{10} is
   reversible for the process $(\eta_{[s,t],\Lambda},\, t\ge s)$ for each fixed
   $s\in\R$.
\end{iten}

\begin{iten}
  \label{i6} For any configuration $\z\in\XX $ and fixed $t\in \R$, the limit
  \begin{equation}
    \label{f6}
    \lim_{s\to-\infty} \eta_{[s,t],\Lambda}(i) :=\eta_{t,\Lambda}(i)
  \end{equation}
  exists with probability one and \emph{does not depend on $\z$}. The process
  $(\eta_{t,\Lambda},\,t\in\R)$ is stationary with time-marginal
  $\mu_\Lambda(\cdot|\y_{\Lambda^c})$, which is the unique invariant measure for
  the process.
\end{iten}

\begin{iten}
  \label{i7} For any configuration $\y\in\XX $ and fixed $t\in \R$, the limit
  \begin{equation}
    \label{f7p}
    \lim_{\Lambda\nearrow\Z} \eta_{t,\Lambda}(i) =\eta_{t}(i)
  \end{equation}
  holds with probability one and \emph{does not depend on $\y$}.
\end{iten}

Theorem
\ref{tp6} follows from \ref{i3} to \ref{i7}..

\section{Finite volume}
\label{S4}
In this section we prove \ref{i5}, \ref{i51} and \ref{i6} for the finite
process $\eta_{t,\Lambda}$.

\paragraph{Generator} {\it Proof of \ref{i5}.}
{}From \eqref{f31} it follows that
\begin{align}
  \label{f311} 
&\eta_{[s,u+h],\Lambda}(i) \;=\; \one\{ |\cN([u,u+h],\{i\})|=0
  \}\,\eta_{[s,u],\Lambda}(i) + \one\{
  |\cN([u,u+h],\{i\})|=1 \}\\
  &\qquad\times\Bigl(\alpha\sum_{j\in\Lambda} p(i,j)\eta_{[s,u],\Lambda}(j) +
  \alpha\sum_{j\in\Lambda^{c}} p(i,j)\y(j)+(1-\alpha)\dd(i)+\xi(\tau,i)\Bigr) 
  + \hbox{other terms}\nonumber 
\end{align}
where $\tau\in (u,u+h)$, $\{ |\cN([u,u+h],\{i\})|=k \}$ is the event ``there are
exactly $k$ Poisson epochs in $[u,u+h]\times\{i\}$'' and the ``other terms'' are
related to the presence of more than one Poisson epoch in $ U(i)\times[u,u+h]$
which has probability of order $h^2$, where $U(i)$ is the cube centered at $i$
with side $R$. The independence properties of the Poisson process and
\eqref{f311} show that the process is Markovian with generator \eqref{p2f}.
\hfill\qed

\paragraph{Invariant measure} {\it Proof of \ref{i51}.}
  Let $p(0,0)=0,\quad 0<h<1$ and $J_{i,j}=(1-h)p(i,j)$ and $h =1-\alpha$. We
  prove that $\mu_{\Lambda}$ defined by \eqref{10} on $\R^\Lambda$ is reversible
  for each one of the generators $L_k$ defined in \eqref{p2f}, that is, for
  $f,g: \R^\Lambda\to\R$,
\begin{equation}\label{35.1}
\mu_{\Lambda}(gL_{k}f)=\mu_{\Lambda}(fL_{k}g).
\end{equation}
For a fixed configuration $\z$ and site $k$ let
$k\in\Lambda$, let
\begin{equation}
  \label{rf10}
 \ovl{\z}(k)=\sum\limits_{i:\:i\neq
k}J_{k,i}\z(i)+h\dd(k)
\end{equation}
Let $\z=\x_{\Lambda} \y_{\Lambda^{c}}$.  A direct calculation yields
\begin{eqnarray}
\label{45.1}
R_{k} &:=&\sum_{i:\:i\neq
k}J_{k,i}(\z(i)-\z(k))^{2}+h(\z(k)-\dd(k))^{2} \\
&=&(\z(k)-\ovl{\z}(k))^{2} +  \sum_{i:\:i\neq k}
          J_{k,i}(\z(i)-\ovl{\z}(k))^{2}+h(\ovl{\z}(k)-\dd(k))^{2}.
          \label{45.1}
\end{eqnarray}
It follows that the conditional law at coordinate $k$ given the heigths at the
other coordinates is a Gaussian with mean $\ovl{\z}(k)$: 
\begin{equation}
  \label{rf9}
  \mu_{\{k\}}(\ed z|\z_{\{k\}^c}) = \frac1{\sqrt {2\pi}}e^{-(z-\ovl{\z}(k))^{2}}
  \ed z
\end{equation}
so that the updating at site $k$ is done with the conditional distribution given
the heights at the other sites. This implies reversibility; to show it we use
the notations
\begin{align*}
E^{k}(\z)&=\exp\left\{-\sum_{i,j:\:i,j\neq k}
J_{i,j}(\z(i)-\z(j))^{2}-\sum_{i:\:i\neq k}h(\z(i)-\dd(i))^{2}\right\} \\
E_{k}(\z,\z(k))&=\exp\left\{-\sum_{i:\:i\neq
k}J_{k,i}(\z(i)-\z(k))^{2}-h(\z(k)-\dd(k))^{2}\right\}
\end{align*}
We represent the left-hand side in (\ref{35.1}) as
$$\mu_{\Lambda}(gL_{k}f)=S^{1}_{k}-S^{2}_{k}, $$
where
\begin{align*}
S^{1}_{k} &:=\frac{1}{Z_{\Lambda}}\int
\prod_{i\in\Lambda}\ed\z(i)
E^{k}(\z)E_{k}(\z,\z(k)) \int_{\mathbb{R}}\ed
           xe^{-x^{2}}g(\z)f\left((1-h)P_{k}\z+
           h\dd(k)+\e_{k}x\right) \\
S^{2}_{k}&:=\frac{\sqrt{\pi}}{Z_{\Lambda}}\int
\prod_{i\in\Lambda}\ed\z(i) E^{k}(\z)E_{k}(\z,\z(k))g(\z)f(\z)
\end{align*}
Using (\ref{45.1}), we obtain
\begin{equation*}
\begin{aligned}
S^{1}_{k}= \frac{1}{Z_{\Lambda}}\int
   \prod_{i\in\Lambda}\ed\z(i)
E^{k}(\z)E_{k}(\z,\ovl{\z}(k))
       \exp\left\{-(\z(k)-\ovl{\z}(k))^{2}\right\}\times \\
       \int_{\mathbb{R}}\ed xe^{-x^{2}}g(\z)f\left((1-h)P_{k}\z+
       h\dd(k)+\e_{k}x\right)
\end{aligned}
\end{equation*}
Next we change the variables as follows
\begin{equation*}
\vb(i)=  \left\{
        \begin{array}{ll}
             \ovl{\z}(k)+x, & \hbox{\rm if\ } i = k; \\
             \vb(i) =\z(i), & \hbox{\rm otherwise}; \\
        \end{array}
        \right.
\end{equation*}
and
\begin{equation*}
y=\z(k)-\ovl{\z}(k).
\end{equation*}

Observe that $\ovl{\vb}(k)=\ovl{\z}(k)$. Therefore
$\z(i)-\ovl{\z}(k)=\vb(i)-\ovl{\vb}(k)$ for any $i\neq k$.
Using first the above substitution and then relation (\ref{45.1})
again yields
\begin{eqnarray*}
S^{1}_{k}
&=& \frac{1}{Z_{\Lambda}}
          \int\prod_{i\in\Lambda}\ed\vb(i)
E^{k}(\vb)E_{k}(\vb,\ovl{\vb}(k))\exp\left\{-(\vb(k)-\ovl{\vb}(k))^{2}
        \right\} \\
& & \quad \quad \times \int_{\mathbb{R}}\ed ye^{-y^{2}}
    g((\ovl{\vb}(k)+y) \vb_{\Lambda \setminus k})f\left(\vb\right)\\
&=& \frac{1}{Z_{\Lambda}}\int \prod_{i\in\Lambda}\ed\vb(i)
E^{k}(\vb)E_{k}(\vb,\vb(k))\int_{\mathbb{R}}\ed
     ye^{-y^{2}}g\left((1-h)P_{k}\vb+ h\dd(k)+\e_{k}y\right)f(\vb).
\end{eqnarray*}
We thus obtain
$\mu_{\Lambda}(fL_{k}g)=S^{1}_{k}-S^{2}_{k} = \mu_{\Lambda}(g L_{k} f)$.
\hfill\qed

\paragraph{Limiting Distributions}
{\it Proof of \ref{i6}.}
We consider separately the four terms in \eqref{f31}. Define
\begin{align}
  \label{f31a}
  \aA_{[s,t],\Lambda}(i) &:=
  \sum_{(j,\tau)\in \cN([s,t],\Lambda)}b_{[\tau,t],\Lambda}(i,j)\xi(j,\tau)\\
  \label{f31b}
  \aB_{[s,t],\Lambda}(i) &:=
  \sum_{(j,\tau)\in \cN([s,t],\Lambda)} a_{[\tau,t],\Lambda}(i,j^*)\dd(j)\\
  \label{f31c}
  \aC_{[s,t],\Lambda}(i) &:= \sum_{(j,\tau)\in \cN([s,t],\Lambda)}
  \sum_{k\in \Lambda^c}a_{[\tau,t],\Lambda}(i,k)\,\y(k)\\
  \label{f31d}
  \aD_{[s,t],\Lambda}(i) &:= \sum_{j\in\Z}b_{[s,t],\Lambda}(i,j) \z(j)
\end{align}
Since the variables $\xi$ are independent of the variables $b$,
$(\aA_{[t-u,t],\Lambda}(i),\,u\ge 0)$ is a martingale for the sigma field
$\cF_u$ generated by $\cN([t-u,t],\Lambda)$ and the associated Gaussian
variables~$\xi$. The variance at time $u$ is given by
\begin{equation}
  \label{f9}
  \E (\aA_{[t-u,t],\Lambda}(i))^2
  = \E \sum_{(j,\tau)\in \cN([t-u,t],\Lambda)} (b_{[\tau,t],\Lambda}(i,j))^2
\end{equation}
by conditioning on the Poisson epochs and deleting the cross terms by
independence of the different  $\xi$.  The last expression is bounded by
\begin{equation}
  \label{f10}
  \le \E \sum_{(j,\tau)\in \cN([t-u,t],\Lambda)} b_{[\tau,t],\Lambda}(i,j) =
  \sum_j \int_0^u e^{-(1-\alpha)r} Q^r_\Lambda (i,j)\ed r \le \frac{1}{1-\alpha}
\end{equation}
uniformly in $\Lambda$, where $Q^r_\Lambda (i,j)$ is the probability for a
continuous time random walk with rates $p$ absorbed at $\Lambda^c$ starting at
$i$ to be at $j$ by time $r$. The factor $e^{-(1-\alpha)r}$ corresponds to the
probability that the walk is not killed in the time interval $[t,t-r]$. The
Martingale convergence theorem implies that
\begin{equation}
  \label{f11}
  \lim_{u\to\infty} \aA_{[t-u,t],\Lambda}(i) = \aA_{t,\Lambda} (i) \quad\hbox{ a.s.}
\end{equation}
The process $(\aA_{t,\Lambda},\,t\in\R)$ is Markov, stationary with generator
\begin{equation}
  \label{35a} \sum_k L^A_{k}f(\x_{\Lambda} \y_{\Lambda^{c}})=\sum_k \int G(\ed x)
  \Big[f\big(\alpha P_{k}(\x_{\Lambda} \y_{\Lambda^{c}}+\e_k x) \big) - f(\x_{\Lambda}
  \y_{\Lambda^{c}})\Big],
\end{equation}
and $\E (\aA_{t,\Lambda} (i))^2\le (1-\alpha)^{-1}$.

We show now that the limit in \eqref{f31b} is given by
\begin{equation}
  \label{f13}
  \lim_{u\to\infty} \aB_{[t-u,t],\Lambda} (i) = \sum_{j\in\Lambda}
  \kK_\Lambda(i,j) \dd(j):= \m_\Lambda(i)
\end{equation}
where $K_\Lambda(i,j)$ is the probability the random walk to be killed at site 
$j$, if $j\in \Lambda$, or to be absorbed, if $j\in \Lambda^c$.
As in \eqref{hh1}, $\m_\Lambda$ satisfies
\begin{equation}
  \label{hh2}
  \m_\Lambda(i)\,=\, \sum_{k\in\Lambda}  \alpha p(i,k)\m_\Lambda(k) + (1-\alpha) \dd(i) .
\end{equation}
For $v>u$, the following ``Markov property'' holds
\begin{equation}
  \label{f12a}
  \aB_{[t-v,t],\Lambda}(i):= \aB_{[t-u,t],\Lambda}(i)
  + \sum_{k\in \Lambda} b_{[t-u,t],\Lambda}(i,k)\aB_{[t-v,t-u],\Lambda}(k)
\end{equation}
Applying \eqref{hh2} to each Poisson epoch $(j,\tau)$ we have
\begin{equation}
  \label{f12b}
  \m_\Lambda(i) = \aB_{[t-u,t],\Lambda}(i)+ \sum_{k\in \Lambda} b_{[t-u,t],\Lambda}(i,k) \m_\Lambda(k)
\end{equation}
which implies that $\m_\Lambda$ is invariant for this dynamics.  For fixed $s$
the process $(\aB_{[s,t],\Lambda},\,t\ge s)$ on $\R^\Lambda$ is Markov with
generator
\begin{equation}
  \label{35b} \sum_{k\in \Lambda} L^B_{k}f(\x_{\Lambda})=\sum_{k\in \Lambda}
  \Big[f\big([\alpha P_{k}(\x_{\Lambda} \oo_{\Lambda^{c}})+(1-\alpha)\e_k\dd(k)]_\Lambda \big) - f(\x_{\Lambda})\Big],
\end{equation}
where the $\oo$ is the ``all zero'' configuration: $\oo(i)\equiv 0$. Finally,
using \eqref{f12b},
\begin{equation}
  \label{f14}
  |\aB_{[t-u,t],\Lambda}(i)- \m_\Lambda(i)| \le \sum_k b_{[t-u,t],\Lambda}(i,k) |\m_\Lambda(k)|
\end{equation}
By \eqref{f10} the right hand side converges to zero as $u\to\infty$ uniformly
in $\Lambda$, proving \eqref{f13}.

A similar argument shows that
\begin{equation}
  \label{f18}
  \lim_{u\to\infty} \aC_{[t-u,t],\Lambda} (i) = \sum_{j\in\Lambda^c}
  \kK_\Lambda(i,j) \y(j):= \rr_\Lambda(i)
\end{equation}

Finally, by \eqref{f10},
\begin{equation}
  \label{f19}
  \lim_{u\to\infty} \aD_{[t-u,t],\Lambda} (i) = 0
\end{equation}

The limits \eqref{f11}, \eqref{f13}, \eqref{f18} and \eqref{f19} show
\eqref{f6}. The resulting limit is
\begin{equation}
  \label{f20}
  \eta_{t,\Lambda} = A_{t,\Lambda} + \m_\Lambda + \rr_\Lambda
\end{equation}
By construction the law of this limit does not depend on $t$. One proves that
the process $(\eta_{t,\Lambda},\,t\in\R)$ is Markov like in \ref{i5}. Since the
limit \eqref{f20} does not depend on the initial configuration $\z$, this shows
that the process $(\eta_{t,\Lambda})$ has a unique invariant measure. Since by
\ref{i6} $\mu_\Lambda(\cdot|\y_{\Lambda^c})$ is invariant for this process, for
each $t\in\R$ the marginal law of $\eta_{t,\Lambda}$ is
$\mu_\Lambda(\cdot|\y_{\Lambda^c})$.
\hfill\qed

\section{Infinite volume and thermodynamic limit}
\label{S5}
\paragraph{Existence of the process} {\it Proof of \ref {i1}.}
The fact that the sums in \eqref{3.1} are finite follows immediately from the
finite range condition on $p$. The proof that the dynamics is the harness
process follows then from an argument similar to the one in
\ref{i5}. \hfill\qed

\paragraph{Limiting stationary processes}
Here we prove \ref{i3} and \ref{i7}. To prove \ref{i3} we consider separately
the three terms in \eqref{3.1}. Define
\begin{align}
 \label{f41a}
\aA_{[s,t]}(i) &:=
  \sum_{(j,\tau)\in \cN([s,t])}b_{[\tau,t]}(i,j)\xi(j,\tau)\\
  \label{f41b}
\aB_{[s,t]}(i)  &:=
  \sum_{(j,\tau)\in \cN([s,t])} a_{[\tau,t]}(i,j^*)\dd(j)\\
 \label{f41d}
 \aD_{[s,t]}(i) &:= \sum_{j\in\Z}b_{[s,t]}(i,j) \z(j)
\end{align}

\begin{lemma}
  \label{f33}
Under the hypothesis of Theorem \ref{tp6},
\begin{equation}
  \label{f24}
  \lim_{u\to\infty} \aA_{[t-u,t]}(i) = \aA_t(i)
\end{equation}
and $(\aA_t,\,t\in\R)$ is a stationary Markov process with generator
\begin{equation}
  \label{45a} \sum_k L^\aA_{k}f(\x)=\sum_k \int G(\ed x)
  \Big[f\big(\alpha P_{k}\x+\e_k x \big) - f(\x)\Big],
\end{equation}
with $\E (\aA_{t} (i))^2\le (1-\alpha)^{-1}$. Furthermore,
\begin{equation}
  \label{f312}
  \lim_{\Lambda\nearrow\Z} \aA_{t,\Lambda}(i) = \aA_t(i)
\end{equation}
\end{lemma}

\begin{proof}
  As in the proof of \eqref{f11}, $(\aA_{[t-u,t]}(i),\,u\ge 0) $ is a martingale
  with uniformly bounded second moments. Hence the limit \eqref{f24} exists.
  Since $b_{ [\tau,t],\Lambda}(i,j)\nearrow b_{[\tau,t]}(i,j)$,
  \begin{equation}
    \label{f37}
    |\aA_t(i)-\aA_{t,\Lambda}(i)|
    \le \sum_{(j,\tau)\in\cN[-\infty,t]}(b_{[\tau,t]}(i,j)-b_{[\tau,t],\Lambda}(i,j)) |\xi_{(j,\tau)}|
  \end{equation}
  The sum is finite because the variance of $\aA_t(i)$ is bounded and it
  converges to zero by monotone convergence.

  The process $(\aA_t,\,t\in\R)$ is stationary by construction. It is Markov
  with generator \eqref{45a} as a consequence of \eqref{f312} and \ref{i5}.
\end{proof}

\begin{lemma}\label{2lp2} If $\dd\in\XX$ then $|\m(i)|<\infty$ for all
  $i\in\Z$ and
  \begin{equation}
    \label{f22}
    \lim_{u\to\infty} \aB_{[t-u,t]}(i) = \m(i)
  \end{equation}
and
\begin{equation}
  \label{f25}
  \lim_{\Lambda\nearrow\Z} \m_\Lambda(i) = \m(i)
\end{equation}
\end{lemma}
\begin{proof}
  Remark that $\kK(i,j)\leq\alpha^{\left[\frac{|i-j|}{R}\right]}$ because at least
  $\left[\frac{|i-j|}{R}\right]+1$ jumps are needed to achieve the point $j$
  from $i$ and during the time $\sigma_{[-\infty,t]}^i$ is not absorbed.
  \begin{equation}
    \label{f21}
    |\m(i)| \le \sum_j \kK(i,j) |\dd(j)| \le \sum_j
    \alpha^{\left[\frac{|i-j|}{R}\right]}  |\dd(j)|<\infty
  \end{equation}
  by the definition of $\XX$. The proof of \eqref{f22} works as the proof of
  \eqref{f13}.

On the other hand $\kK_\Lambda(i,j)\nearrow \kK(i,j)$, hence
\begin{equation}
  \label{f30}
  |\m(i)-\m_\Lambda(i)| \le \sum_j(\kK(i,j)-\kK_\Lambda(i,j)) |\dd(j)| \to 0
\end{equation}
by monotonic convergence.
\end{proof}

\begin{lemma}
\label{exist}
If $\y\in\XX$ then for $\rr_\Lambda$ defined in \eqref{f18},
\begin{equation}
  \label{f26}
  \lim_{\Lambda\nearrow\Z} \rr_\Lambda(i) = 0
\end{equation}
\end{lemma}
\begin{proof}
As in the proof of the previous lemma,
  \begin{equation}
    \label{f27}
    |\rr_\Lambda(i)| \le \sum_{j\in\Lambda^c} \alpha^{\left[\frac{|i-j|}{R}\right]}  |\y(j)|
  \end{equation}
  which converges to zero as $\Lambda\nearrow\Z$ because the sum is finite by
  definition of $\XX$.
\end{proof}

\paragraph{\it Proof of \ref{i3}} Calling $\eta_t=: A_t+\m$, the convergence
follows from \eqref{f24}, \eqref{f22} and \eqref{f19} (which holds uniformly in
$\Lambda$). \hfill\qed

\paragraph{\it Proof of \ref{i7}}
Since
\begin{equation}
\label{f40}
\eta_{t,\Lambda} = \aA_{t,\Lambda} + \m_\Lambda +\rr_\Lambda \quad\hbox{ and
}\quad\eta_{t} = \aA_t + \m \,,
\end{equation}
the convergence follows from \eqref{f312}, \ref{f25} and \ref{f26}.
\hfill\qed

\paragraph{Limiting process satisfies DLR conditions} {\it Proof of \ref{i4}.} 
Here we use the thermodynamic limit. The measures
$\mu_\Lambda(\cdot|\y_{\Lambda^c})$ are compatible with the DLR conditions for
subsets of $\Lambda$. By \eqref{f7p} the measure $\mu$ is the limit as $\Lambda$
increases to $\Z$ of $\mu_\Lambda(\cdot|\y_{\Lambda^c})$, independently of
$\y\in\XX$. Hence, also $\mu$ is compatible with the DLR conditions. \hfill\qed

\paragraph{DLR measures are invariant} {\it Proof of \ref{i4a}.} 
This proof does not use the thermodynamic limit. We first prove that $\tmu$
is invariant.  Let $\eta$ be a configuration sampled with $\tmu$. Let
$\ovl\eta(k)$ be defined in \eqref{rf10}. Let $Z$ be a standard Gaussian
variable independent of $\eta$. Since $\tmu$ satisfies the DLR equations, a
computation like in \eqref{45.1} shows that the variable
\[
\ovl\eta(k) +Z \;\;\hbox{ has law }\;\; \tmu (\x(k)\in\cdot
\,|\,\x_{\{k\}^c}=\eta_{\{k\}^c})
\]
 and the configuration
\[
\theta(k,\eta,Z) = \eta - \delta_{\eta(k)} + \delta_{\ovl\eta(k) +Z }
\]
obtained by substituting the value at $k$ by $\ovl\eta_k +Z$ has law $\tmu$. The
Harness dynamics does the same substitution with the Gaussian noises $\xi$ at
the updating epochs; see \eqref{f7}. This means that after each Poisson epoch,
the distribution of the updated configuration is the same as the one just before
the updating. But since there are infinitely many sites, there is no ``first
Poisson epoch'' to apply the rule and show \ref{i4a} directly. To overcome the
difficulty consider a time $t$ small enough such that $\Z$ is partitioned in
finite sets $\Lambda_\ell$, $\ell=1,2,\dots$ (depending on the Poisson epochs in
the interval $[0,t]$) in such a way that if there is an epoch at time
$s\in[0,t]$ at site $i\in\Lambda_\ell$, then $\{j\in\Z\,:\, ||i-j||\le R\}\subset
\Lambda_\ell$, where we recall $R$ is such that $p(0,j)=0$ if $||j||\ge R$. A
standard percolation argument shows that this is possible.

Enumerate the epochs as follows. Start ordering cronologically the epochs in
$\Lambda_1$; continue with the epochs in $\Lambda_2$ and so on. Under this
ordering, if $n(i,s)$ is the label of epoch $(i,s)$ then $n(i,s)<n(i',s')$ if
$s<s'$ or if $i\in \Lambda_\ell$, $i'=\Lambda_{\ell'}$ and $\ell<\ell'$.  Call
$(i_n,s_n)$ the $n$th Poisson epoch in this labeling and $\xi_n=\xi(i_n,s_n)$
the associated Gaussian variable. Choose $\eta$ with law $\tmu$ and define
inductively $\eta_0=\eta$ and for $n\ge1$,
\[
\eta_n = \theta(i_n,\eta_{n-1},\xi_n)
\]
By the previous considerations $\eta_n$ has law $\tmu$ for all $n$. By the
definition of $\Lambda_\ell$, the updating of sites in $\Lambda_\ell$ in the
time interval $[0,t]$ depends on the initial configuration $\eta$ only through
the values at sites in $\Lambda_\ell$. On the other hand, for any cylinder
function with support on a finite set $\Lambda$, there exists an
$n(\Lambda,\cN)<\infty$ such that $\Lambda\subset
\{i_1,\dots,i_{n(\Lambda,\cN)}\}$ and $\int \tmu(\ed \eta)\E
(f(\eta_{n(\Lambda,\cN)})\,|\,\cN)=\tmu f$ for almost all $\cN$. Hence
\[
\tmu S(t) f 
\;=\; \int \tmu(\ed\eta)\,\E\big(\E (f(\eta_{n(\Lambda,\cN)})\,|\,\cN)\big) 
\;=\; \E\Big(\int \tmu(\ed\eta)\,\E (f(\eta_{n(\Lambda,\cN)})\,|\,\cN)\Big) 
\;=\; \tmu f
\] 
This shows invariance of $\tmu$. 

A similar argument shows reversibility. Going backwards in time, the law
of $\eta_{n-1}(i_n)-\bar\eta_n(i_n)$ is a standard Gaussian
variable independent of $\eta_n$. The construction can be done using
the same Poisson epochs and partition $(\Lambda_\ell)$.\hfill \qed

\section{Zero temperature}
\label{S6}
\begin{proof}[Proof of Theorem \ref{14}]
  We need to show that the configuration $\m$ defined in \reff{15} is compatible
  with the specifications \reff{10} at zero temperature.  Let $\Lambda\subset\Z$
  be a finite volume. We want to show that $\m_\Lambda$ given by \eqref{15} is
  the $\x_\Lambda$ which minimizes
\begin{align}\label{Hham}
H(\x_\Lambda\m_{\Lambda^c})=\alpha\sum_{i,j\in \Lambda}p(i,j)(\x(i)-\x(j))^{2}+
\alpha\sum_{\begin{smallmatrix} i\in \Lambda\\
k\in \Lambda^{c}
\end{smallmatrix}}p(i,k)(\x(i)-\m(k))^{2}\\
+(1-\alpha)\sum_{i\in \Lambda}(\x(i)-\dd(i))^{2},\nonumber
\end{align}
The Hamiltonian \reff{Hham} can be
rewritten as
\begin{align}
\label{4g1}
&H(\x_\Lambda\m_{\Lambda^c})\;=\;\alpha\sum_{i\in \Lambda}\sum_{k\in \Z}
p(i,k)(\x(i)-\m(k))^{2}+
(1-\alpha)\sum_{i\in \Lambda}(\x(i)-\dd(i))^{2}\nonumber\\
&\quad+\,2\alpha\sum_{i,j\in \Lambda}
p(i,j)\Big[-(\x(i)-\m(j))(\x(j)-\m(i))\\
&\quad+\,(\x(i)-\m(j))(\m(j)-\m(i))-(\x(j)-\m(i))(\m(j)-\m(i))+(\m(j)-\m(i))^{2}\Big]\nonumber
\end{align}
Let $r_{i}=\m(i)^{2}-\Big[\alpha\sum_{k\in
\Z} p(i,k)\m(k)^{2}+(1-\alpha)\dd(i)^{2}\Big]$ then
\begin{align*}
  &H(\x_\Lambda\m_{\Lambda^c})\;= \; \sum_{i\in
    \Lambda}(\x(i)-\m(i))^{2}-\sum_{i\in \Lambda}r_{i}\\ 
  &\quad +\,2\alpha\sum_{i,j\in \Lambda}p(i,j)\Big[-(\x(i)-\m(j))(\x(j)-\m(i))+(\x(i)-\m(j))(\m(j)-\m(i))\\
  &\qquad\qquad\qquad\qquad\qquad\qquad -\,(\x(j)-\m(i))(\m(j)-\m(i))+(\m(j)-\m(i))^{2}\Big]
\end{align*}
To find the configuration $\x_\Lambda$ minimizing \eqref{Hham} we take the
derivatives for $i\in\Lambda$:
\begin{align*}
  \frac{\partial H}{\partial \x(i)} = &\;2(\x(i)-\m(i)) \\
  &+\,\alpha\sum_{k\in \Lambda} p(i,k)\Big[-2(\x(k)-\m(i))+2(\m(k)-\m(i))\Big]
\end{align*}
Then we obtain the following system of the equations
\begin{equation}\label{5g1}
\x(i)-\alpha\sum_{k\in  \Lambda}p(i,k)\x(k)=
\m(i)-\alpha\sum_{k\in  \Lambda}p(i,k)\m(k).
\end{equation}
It is clear that $\x(i)=\m(i)$ is a solution of \reff{5g1}. Since $H_\Lambda$ is convex,
$\m$ is the minimum. Since $H_{\Lambda}$ is a second order polynomial the
minimum is unique. This shows that $\m$ belongs to the support of a ground state
measure.

\paragraph{Uniqueness} Let $\tm\in\XX$ be a minimizer of $H$. Then 
$\tm(i)=\alpha\sum_j p(i,j)\tm(j) +(1-\alpha) \dd(j)$ minimizes
$H(\x_{\{i\}}\tm_{\{i\}^c})$ for each $i\in\Z$. So $\tm$ is invariant for the dynamics
\begin{equation}
  \label{35bb} \sum_{k\in \Z} L^B_{k}f(\x)=\sum_{k\in \Z}
  \Big[f\big(\alpha P_{k}(\x)+(1-\alpha)\e_k\dd(k) \big) - f(\x)\Big],
\end{equation}
because this dynamics chooses a site $i$ at random times and substitutes the
height at $i$ with the value minimizing $H(\x_{\{i\}}\tm_{\{i\}^c})$. The
graphical construction of this process with initial measure $\tm$ and the
invariance of $\tm$ gives
\begin{equation}
  \label{f12bb}
  \tm(i) = \aB_{[t-u,t]}(i)+ \sum_{k\in \Z} b_{[t-u,t]}(i,k) \tm(k)
\end{equation}
Subtracting this equation from \eqref{f12b} with $\Lambda=\Z$ we get
\begin{equation}
  \label{f12bp}
  | \m(i)-\tm(i)| \le  \sum_{k\in \Z} b_{[t-u,t]}(i,k)| \m(k)-\tm(k)|
\end{equation}
This is summable because both $\m$ ad $\tm$ belong to $\XX$ and since
$\lim_{u\to\infty}b_{[t-u,t]}(i,k)=0$, then $\m=\tm$.
\end{proof}

\section*{Acknowledgements}
We thank Roberto Fernandez and a careful referee for discussions about
uniqueness of Gibbs states and many other comments that improved the paper.

This paper has been partially supported by Funda\c c\~ao de Amparo \`a Pesquisa
do Estado de S\~ao Paulo (FAPESP) and Conselho Nacional de Desenvolvimento
Cient\'{\i}fico (CNPq), Brazil. E.P. was partially supported by Grant CRDF
RUM1-2693--MO-05.

\parindent0pt
\parskip 0pt

Pablo A.~Ferrari, Beat M.~Niederhauser,\\
IME USP,
Caixa Postal 66281,
05311-970 - S\~{a}o Paulo,
BRAZIL\\
Phones: +55 11 3091 6119, +55 11 3091 6129, Fax: +55 11 3814 4135\\
{\tt pablo@ime.usp.br},
{\tt http://www.ime.usp.br/\~{}pablo}
\vskip 1mm
\vskip 1mm

Eugene A.~Pechersky\\
IITP,
19, Bolshoj Karetny per.,
GSP-4,
Moscow  127994,
RUSSIA\\
{\tt pech@iitp.ru}
\end{document}